\documentclass[11pt]{article}
\usepackage{cite}
\usepackage{amsmath}
\usepackage{amsfonts}
\usepackage{amssymb}
\usepackage{stmaryrd}
\usepackage{graphicx}
\usepackage{subfigure}
\usepackage{caption}
\usepackage{color}

\numberwithin{equation}{section}

\newcommand{\sech}{{\mathrm{sech}} }
\newcommand{\p}{\partial}
\newcommand{\og}{\omega}

\newcommand{\fl}[2]{\frac{#1}{#2}}
\newcommand{\be}{\begin{equation}}
\newcommand{\ee}{\end{equation}}
\newcommand{\nn}{\nonumber}

\newtheorem{theorem}{Theorem}[section]

\newtheorem{lemma}[theorem]{Lemma}

\begin{document}

\title{A Lawson-type exponential integrator for the Korteweg--de Vries equation}

\author{Alexander Ostermann\thanks{Department of Mathematics,
University of Innsbruck, 6020 Innsbruck, Austria\newline ({\tt alexander.ostermann@uibk.ac.at})}\,,
Chunmei Su\thanks{Department of Mathematics,
University of Innsbruck, 6020 Innsbruck, Austria ({\tt sucm13@163.com})}}
\date{}
\maketitle

\begin{abstract}
We propose an explicit numerical method for the periodic Korteweg--de Vries equation. Our method is based on a Lawson-type exponential integrator for time integration and the Rusanov scheme for Burgers' nonlinearity. We prove first-order convergence in both space and time under a mild Courant--Friedrichs--Lewy condition $\tau=O(h)$, where $\tau$ and $h$ represent the time step and mesh size, respectively, for solutions in the Sobolev space $H^3((-\pi, \pi))$. Numerical examples illustrating our convergence result are given.
\end{abstract}


\pagestyle{myheadings}\thispagestyle{plain}

\section{Introduction}
Consider the Korteweg--de Vries (KdV) equation
\be\label{KdV}
\begin{aligned}
& u_t+u_{xxx}+u u_x=0, \quad x \in \Omega=(-\pi, \pi),\quad t>0,\\
&u(x,0)=u_0(x),\quad x \in \overline{\Omega},
\end{aligned}
\ee
where we impose periodic boundary conditions for practical implementation. The KdV \,equation is a generic model for the study of weakly nonlinear long waves. It describes the propagation of shallow water waves in a channel \cite{Kordeweg1895} and is widely applied
in science and engineering, such as in plasma physics where it gives rise to ion acoustic solitons \cite{Das1998} and in geophysical fluid dynamics where it describes long waves in shallow seas and deep oce\-ans \cite{Osborne1995, Ostrovsky1989}. The KdV equation is also relevant for studying the interaction between nonlinearity and dispersion.

For the well-posedness of the periodic KdV equation, we refer to \cite{Bourgain1993, Tao2003, Gubinelli2012}. It was shown in \cite{Tao2003} that the equation is globally well-posed for initial data in $H^s(\Omega)$ with $s\ge-1/2$. For its numerical solution, various methods have been proposed and analyzed in the literature, such as finite difference methods (FDM) \cite{Vliegenthart1971, Taha1984, Helal2007, Holden2014}, finite element methods \cite{Aksan2006, Winther1980, dutta2015, Arnold1982}, Fourier spectral methods \cite{Chan, Maday, Rashid, Rashid2006, Rashid2007}, splitting methods \cite{Holden1999, Klein2008, Holden2011} and Petrov--Galerkin methods for the KdV equation with nonperiodic boundary condition \cite{Ma2000legendre, Ma2001optimal, Shen2003}. Numerical methods for the Kadomtsev--Petviashvili equation, which is a two-dimensional generalization of the KdV equation, were considered in \cite{Lukas2017, Klein2011fourth}.

For finite difference methods, linear stability has been analyzed in \cite{Goda1975, Taha1984, Vliegenthart1971}. The explicit leap-frog scheme \cite{Taha1984} and the Lax--Friedrichs scheme \cite{Vliegenthart1971} require both the rather severe stability condition $\tau=O(h^3)$, where $\tau$ and $h$ represent the discretization parameters in time and space, respectively. To weaken the stability restriction, some implicit FDM were proposed in \cite{Goda1975, Taha1984}. Recently, the Lax--Friedrichs scheme with an implicit dispersion was proved to converge uniformly to the solution of the KdV equation for initial data in $H^3$ under the stability condition $\tau=O(h^{3/2})$ for both the decaying case on the full line and the periodic case \cite{Holden2014}. However, no convergence rate was obtained. Very recently, for the $\theta$-right winded FDM, which applies the Rusanov scheme for the hyperbolic flux term and a 4-point $\theta$-scheme for the dispersive term, first-order convergence in space was proved under a hyperbolic Courant--Friedrichs--Lewy (CFL) condition $\tau=O(h)$ for $\theta\ge\fl{1}{2}$ and under an \emph{Airy} Courant--Friedrichs--Lewy condition $\tau=O(h^3)$ for $\theta<\fl{1}{2}$, for solutions in $H^6(\mathbb{R})$ \cite{Courtes2017}.

On the other hand, the numerical approximation by Fourier spectral/pseudospectral methods has been studied by many authors \cite{Maday, Ma1986}. Maday and Quarteroni \cite{Maday} showed that for solutions in $H^r$, the error of the Fourier spectral method is of order $O(h^{r-1})$ in the $L^2$ norm while the error of the pseudospectral method is of order $O(h^{r-2})$ in the $H^1$ norm. The corresponding $L^2$ estimate for the Fourier pseudospectral method was established in \cite{Ma1986} with the aid of artificial viscosity to avoid the nonlinear instability caused by the aliasing error. More specifically, first-order convergence in time was shown in \cite{Ma1986} for the fully discrete pseudospectral method under the stability condition $\tau=O(h^3)$ for explicit and $\tau=O(h^2)$ for implicit discretization of the nonlinear term, respectively. For the rigorous analysis of splitting methods, we refer to \cite{Holden2011, Holden2013}.

Nowadays, exponential time integration methods are widely applied for parabolic and hyperbolic  problems \cite{Acta, HS2017, Bao2013}. In particular, a distinguished exponential-type integrator was derived for the KdV equation \cite{HS2017} by using a ``twisting" technique. For this integrator, first-order convergence in time was proved without any CFL condition required. However, the success of this scheme strongly depends on the particular form of the equation. The resulting key relation $k_1^3+k_2^3-(k_1+k_2)^3=-3(k_1+k_2)k_1k_2$ in Fourier space allows one to integrate the stiff part involving $\partial_x^3$ exactly without loss of regularity. Such an integrator, however, can hardly be extended to more general equations, e.g., the fifth-order KdV equation, without additional regularity assumptions. Furthermore, the spatial error was not considered in \cite{HS2017}.

In the present paper, we propose a Fourier pseudospectral method based on a classical Lawson-type exponential integrator, which integrates the linear part exactly, and the Rusanov scheme for Burgers' nonlinearity with an added artificial viscosity. The method is explicit, implemented with FFT and efficient in practical computation. First-order convergence in both space and time is shown under a mild CFL condition $\tau=O(h)$. Moreover, the method can be easily extended to other dispersive equations with Burgers' nonlinearity.

The rest of this paper is organized as follows. In Section 2, we present the necessary notation, the numerical scheme and the main convergence result. Section 3 is devoted to the details of the error analysis.
Numerical results are reported in Section 4 to illustrate our error bounds.

Throughout the paper, $C$ represents a generic constant, which is independent of the discretization parameters and the exact solution $u$.

\section{The Fourier pseudospectral method}
We adopt the standard Sobolev spaces and denote by $\|\cdot\|$ and $(\cdot, \cdot)$ the norm and inner product in $L^2(\Omega)$, respectively. For $m\in \mathbb{N}$, we denote by $H_{\rm p}^m(\Omega)$ the $H^m$ functions on the one-dimensional torus $\Omega=(-\pi,\pi)$. In particular, these functions have derivatives up to order $m-1$ that are all $2\pi$-periodic. The space is equipped with the standard norm $\|\cdot\|_m$ and semi-norm $|\cdot|_m$.

Let $\tau=\Delta t>0$ be the time step size and denote the temporal grid points by
$t_k:=k\tau$ for $k=0,1,2,\ldots$. Given a mesh size $h:=2\pi/(2N+1)$ with $N$ being a positive integer, let
$$x_j:=-\pi+jh,\qquad j=0,1,\ldots, 2N,$$
be the spatial grid points in $[-\pi, \pi)$.
Denote
\begin{align*}
&X_N:=\mathrm{span}\left\{e^{ikx}:  |k|\le N\right\},\qquad\qquad\qquad\,
\widetilde{X}_N:=\left\{v=\sum\limits_{k=-N}^{N}v_ke^{ikx}\in\mathbb{R}\right\}\subseteq X_N,\\
&Y_N:=\left\{v=(v_0, v_1, \ldots, v_{_{2N}})\in \mathbb{C}^{2N+1}\right\},\qquad
\widetilde{Y}_N=Y_N\cap \mathbb{R}^{2N+1}.
\end{align*}
For any $u, v\in C(\Omega)$, define the following discrete inner product and norm by
\[\langle u, v\rangle_{_N}=h\sum\limits_{j=0}^{2N} u(x_j)\overline{v(x_j)},\qquad \interleave  u\interleave_{_N}=\langle u, u\rangle_{_N}^{1/2}.\]
For a periodic function $v(x)$ and a vector $v\in Y_N$, let $P_N: L^2(\Omega)\rightarrow X_N$ be the standard orthogonal projection operator, and $I_N: C(\Omega)\rightarrow X_N$ or $I_N: Y_N\rightarrow X_N$ be the interpolation operator \cite{Shen}, i.e.,
\begin{align*}
&(P_N v, \varphi)=(v, \varphi), \quad \mathrm{for}\ \ \mathrm{all}\quad \varphi\in X_N; \\
&(I_Nv)(x_j)=v(x_j), \quad \mathrm{or}\quad (I_N v)(x_j)=v_j,\quad j=0,\ldots,2N.
\end{align*}
More specifically, $P_N v$ and $I_N v$ can be written as
\[(P_N v)(x)=\sum\limits_{l=-N}^{N} \widehat{v}_l e^{ilx}, \quad (I_N v)(x)=\sum\limits_{l=-N}^{N} \widetilde{v}_l e^{ilx},\]
where $\widehat{v}_l$ and $\widetilde{v}_l$ are the Fourier and discrete Fourier coefficients, respectively, defined as
$$\widehat{v}_l=\fl{1}{2\pi}\int_0^{2\pi}v(x)e^{-ilx}dx,\qquad \widetilde{v}_l=\fl{1}{2N+1}\sum\limits_{j=0}^{2N} v_je^{-ilx_j},\quad l=-N,\dots, N.$$
It was proved in \cite{Shen} that for any $u, v\in C(\Omega)$,
\be\label{inp}
\langle u, v\rangle_{_N}=(I_N u, I_N v),\qquad \interleave  u\interleave_{_N}=\|I_N u\|.
\ee

The semi-discrete pseudospectral method for \eqref{KdV} consists in finding $u_{_N}$ in $\widetilde{X}_N$ such that
\be\label{Ikdv}
\begin{aligned}
& \p_t u_{_N}(x,t)+\p_x^3u_{_N}(x,t)+\fl{1}{2}I_N\left((u_{_N}^2)_x\right)(x,t)=0, \quad x \in \Omega=(-\pi, \pi),\quad t>0,\\
&u_{_N}(x,0)=I_N(u_0)(x),\quad x \in \overline{\Omega}.
\end{aligned}
\ee
Thus, by Duhamel's formula, we have
\[u_{_N}(t_n+\tau)=e^{-\tau\p_x^3}u_{_N}(t_n)-\fl{1}{2}\int_0^\tau e^{-(\tau-s)\p_x^3}I_N\left((u_{_N}^2)_x(t_n+s)\right) ds.\]
By applying the approximation $u_{_N}(t_n+s)\approx u_{_N}(t_n)$ and the first-order Lawson method \cite{Lawson1967, Hochbruck2005}, we get a first-order approximation as
\be\label{sch0}
u_{_N}(t_n+\tau)\approx e^{-\tau\p_x^3}u_{_N}(t_n)-\fl{\tau}{2} e^{-\tau \p_x^3}I_N\left((u_{_N}^2)_x(t_n)\right).
\ee
To ensure the stability, we apply the Rusanov scheme \cite{2009numerical, Courtes2017} for Burgers' nonlinearity, which consists of a centered hyperbolic flux and an added artificial viscosity. The scheme then reads as
\be\label{sch}
\begin{split}
&u_{_N}^{n+1}=e^{-\tau\p_x^3}u_{_N}^n-\fl{\tau}{2} e^{-\tau \p_x^3}I_N\delta_x^0 \left((u_{_N}^n)^2\right)+\fl{c\tau h}{2}e^{-\tau \p_x^3}\delta_x^2u_{_N}^n,\quad n\ge 0,\\
&u_{_N}^0=I_N(u_0),
\end{split}
\ee
where the constant $c$ is the so-called Rusanov coefficient, which has to satisfy a certain condition (cf. \eqref{c0}). Moreover, we have used the notation
\[\delta_x^0v(x)=\fl{v(x+h)-v(x-h)}{2h},\quad \delta_x^2 v(x)=\fl{v(x+h)-2v(x)+v(x-h)}{h^2},\]
where $v(x)=v(x\pm 2\pi)$. Similarly, for a vector $v\in Y_N$, define the standard finite difference operators as
\[\delta_x^0v_j=\fl{v_{j+1}-v_{j-1}}{2h},\quad \delta_x^2v_j=\fl{v_{j+1}-2v_j+v_{j-1}}{h^2},\quad \delta_x^+v_j=\fl{v_{j+1}-v_j}{h},\quad j=0, 1,\dots,2N,\]
with $v_{j\pm(2N+1)}=v_j$ when necessary.

We are now in the position to present the main result of the paper.
\begin{theorem}\label{main}
Assume that the solution of \eqref{KdV} satisfies $u\in C(0,T; H_{\rm p}^3(\Omega))$ and let $c_0>0$ be given by condition \eqref{c0}. Then, for $c>c_0$, there exists $h_0>0$ such that for all $h\le h_0$ and $\tau\le h/c$, the error of scheme \eqref{sch} satisfies
\be\label{error}
\|u_{_N}^n-u(t_n)\|\le M(\tau+h),\qquad n\tau\le T.
\ee
Here both of the constants $M$ and $h_0$ depend on $T$, $c$ and $\|u\|_{L^\infty(0,T; H_{\rm p}^3(\Omega))}$ (cf. \eqref{Md} and \eqref{h0}).
\end{theorem}
\section{Error estimate}
The purpose of this section is to prove Theorem \ref{main}.
\subsection{Some lemmas}
We recall three lemmas from the literature and then prove an additional lemma. All the lemmas are used in the proof of Theorem \ref{main}.
\begin{lemma}\hspace{-1mm}\cite{Shen}.\label{proint}
For any $u\in H_{\rm p}^m(\Omega)$ and $0\le \mu\le m$,
\be\label{pro}
\|P_N u-u\|_\mu \le C h^{m-\mu} |u|_m,\quad \|P_N u\|_m\le C\|u\|_m.
\ee
In addition, if $m>\fl{1}{2}$, then
\be\label{int}
\|I_N u-u\|_\mu \le C h^{m-\mu} |u|_m,\quad \|I_N u\|_m\le C\|u\|_m.
\ee
\end{lemma}

\begin{lemma}[Nikolski's Inequality]\hspace{-1mm}\cite{Shen}.
For any $u\in X_N$ and $1\le p\le q\le \infty$,
\[\|u\|_{L^q} \le \left(\fl{Np_0+1}{2\pi}\right)^{\fl{1}{p}-\fl{1}{q}} \|u\|_{L^p},\]
where $p_0$ is the smallest even integer $\ge p$.
In particular, for $p=2$ and $q=\infty$, we have
\be\label{NI}
\|u\|_{L^\infty}\le h^{-1/2}\|u\|.
\ee
\end{lemma}

\begin{lemma}[Bernstein's Inequality]\hspace{-1mm}\cite{Shen}.
For any $u\in X_N$ and $0\le \mu\le m$,
\be\label{BI}
\|u\|_m \le C h^{\mu-m} \|u\|_{\mu}.
\ee
\end{lemma}

\begin{lemma}\label{lema4}
For $a=(a_0, a_1, \ldots, a_{_{2N}})$, $b=(b_0, b_1, \ldots, b_{_{2N}})\in \widetilde{Y}_N$, we have
\be
\langle a, \delta_x^2b\rangle_{_N}=-\langle \delta_x^+a, \delta_x^+b\rangle_{_N},\label{intp}
\ee
and
\begin{align}
&\sum\limits_{j=0}^{2N}\left(\delta_x^2a_j\right)^2=\frac{4}{h^2}\sum\limits_{j=0}^{2N}
\left[\left(\delta_x^+a_j\right)^2-\left(\delta_x^0a_j\right)^2\right],\label{eq0}\\
&\sum\limits_{j=0}^{2N}a_ja_{j+1}\delta_x^+a_j=-\fl{h^2}{3}\sum\limits_{j=0}^{2N}(\delta_x^+a_j)^3,\qquad\qquad
\sum\limits_{j=0}^{2N}a_{j-1}a_{j+1}\delta_x^0a_j=-\fl{4h^2}{3}\sum\limits_{j=0}^{2N}(\delta_x^0a_j)^3,\label{eq1}\\
&\sum\limits_{j=0}^{2N}\delta_x^2a_j\delta_x^0(ab)_j=-\fl{1}{h^2}\sum\limits_{j=0}^{2N} a_ja_{j+1}\delta_x^+b_j+\fl{1}{h^2}\sum\limits_{j=0}^{2N}a_{j-1}a_{j+1}\delta_x^0b_j.\label{eq2}
\end{align}
\end{lemma}
\emph{Proof.} The identity \eqref{intp} is the discrete version of the integration by parts formula:
\begin{align*}
\langle a, \delta_x^2 b\rangle_{_N}&=\fl{1}{h}\sum\limits_{j=0}^{2N}a_j(b_{j+1}-2b_{j}+b_{j-1})=\fl{1}{h}\sum\limits_{j=0}^{2N}a_j(b_{j+1}-b_{j})-
\fl{1}{h}\sum\limits_{j=0}^{2N}a_j(b_{j}-b_{j-1})\\
&=\sum\limits_{j=0}^{2N}a_j\delta_x^+b_j-
\sum\limits_{j=0}^{2N}a_{j+1}\delta_x^+b_j
=-h\sum\limits_{j=0}^{2N}(\delta_x^+a_j)(\delta_x^+b_j)=-\langle \delta_x^+a, \delta_x^+b\rangle_{_N}.
\end{align*}
The equalities \eqref{eq0}-\eqref{eq2} were established in \cite{Courtes2017} for infinite sequences. By applying the same arguments, we can get \eqref{eq0}-\eqref{eq2} for periodic sequences here. We refer to \cite{Courtes2017} for details.
\hfill $\square$ \bigskip
\subsection{Local error analysis}
We introduce the local truncation error $\xi^{n+1}$ as defect
 \be\label{local}
 \xi^{n+1}=u(t_{n+1})-e^{-\tau\p_x^3}u(t_n)+\fl{\tau}{2}e^{-\tau\p_x^3} \left[\delta_x^0(u(t_n)^2)-ch\delta_x^2u(t_n)\right],\quad n\ge 0.
 \ee
 The local error can be bounded as follows.
 \begin{lemma}\label{localer}
For $u\in C(0,T; H_{\rm p}^3(\Omega))$, we have
 \[\|\xi^{n+1}\|\le M_3\tau^2+M_2\tau h,\]
 where $M_3$ and $M_2$ depend on $\|u\|_{L^\infty(0,T;H^3(\Omega))}$ and $\|u\|_{L^\infty(0,T;H^2(\Omega))}$, respectively.
 \end{lemma}
 \emph{Proof.} We first recall
 \[u(t_{n+1})=e^{-\tau\p_x^3} u(t_n)-\fl{1}{2}\int_0^\tau e^{-(\tau-s)\p_x^3}(u^2)_x(t_n+s)ds\]
 and that $e^{t\p_x^3}$ is a linear isometry for all $t\in\mathbb{R}$. This yields that
 \begin{align*}
\|\xi^{n+1}\|&\le\fl{1}{2}\Big\|\int_0^\tau \left[ e^{-(\tau-s)\p_x^3}(u^2)_x(t_n+s)-e^{-\tau\p_x^3}(u^2)_x(t_n)\right]ds\Big\|\\
&\quad+\fl{\tau}{2}
\left\|e^{-\tau\p_x^3}[(u^2)_x(t_n)-\delta_x^0(u(t_n)^2)]\right\|+\fl{c\tau h}{2}\left\|e^{-\tau\p_x^3}\delta_x^2u(t_n)\right\| \\
&=\Big\|\int_0^\tau(\tau-s)\p_s\big[e^{-(\tau-s)\p_x^3}(uu_x)(t_n+s)\big]ds\Big\|\\
&\quad+\fl{\tau}{2}
\left\|(u^2)_x(t_n)-\delta_x^0(u(t_n)^2)\right\|+\fl{c\tau h}{2}\left\|\delta_x^2u(t_n)\right\|\\
&=:I_1+I_2+I_3.
\end{align*}
For the first part we get
\begin{align*}
I_1&=\Big\|\int_0^\tau(\tau-s)e^{-(\tau-s)\p_x^3}\left[\p_x^3(uu_x)+\p_s(uu_x)\right](t_n+s)ds\Big\|\\
&=\Big\|\int_0^\tau(\tau-s)e^{-(\tau-s)\p_x^3}\left[3u_x\p_x^3u+3(\p_x^2u)^2-u^2\p_x^2u-2uu_x^2\right](t_n+s)ds\Big\|\\
&\le C\tau^2 \sup\limits_{0\le t\le T} \|u\|_3\left[\|u_x\|_{L^\infty}+\|\p_x^2u\|_{L^\infty}+\|u\|_{L^\infty}^2+\|u\|_{L^\infty}
\|u_x\|_{L^\infty}\right]\\
&\le C\tau^2 \sup\limits_{0\le t\le T} (\|u\|_3^2+\|u\|_3^3),
\end{align*}
where we employed equation \eqref{KdV} and the Sobolev imbedding theorem $H^3(\Omega)\hookrightarrow W^{2,\infty}(\Omega)$.
Further, using Taylor expansion and H\"{o}lder's inequality, we have
\begin{align*}
I_2&=\fl{\tau}{4h}\Big\|\int_0^h(h-y)
\left[\p_x^2(u^2)(\cdot+y,t_n)-\p_x^2(u^2)(\cdot-y,t_n)\right]dy\Big\|\\
&\le \fl{\tau h^{1/2}}{4}\Big(\int_{-\pi}^{\pi}\int_0^h \left[\p_x^2(u^2)(\cdot+y,t_n)-\p_x^2(u^2)(\cdot-y,t_n)\right]^2dy dx\Big)^{1/2}\\
&\le \tau h\|\p_x^2(u^2)(t_n)\|\le 2\tau h\left(\|u_x(t_n)\|_{L^\infty}\|u(t_n)\|_1+\|u(t_n)\|_{L^\infty} \|u(t_n)\|_2\right)\\
&\le C\tau h \|u(t_n)\|^2_2\le C\tau h \sup\limits_{0\le t\le T} \|u\|_2^2.
\end{align*}
A similar calculation shows that
\begin{align*}
I_3&=\fl{c\tau}{2h}\Big\|\int_0^h(h-y)
\left[\p_x^2u(\cdot+y,t_n)+\p_x^2u(\cdot-y,t_n)\right]dy\Big\|\\
&\le c\tau h^{1/2}\Big(\int_{-\pi}^{\pi}\int_0^h \left[\p_x^2u(\cdot+y,t_n)+\p_x^2u(\cdot-y,t_n)\right]^2dy dx\Big)^{1/2}\\
&\le 2c\tau h \sup\limits_{0\le t\le T} \|u\|_2,
\end{align*}
which completes the proof.
\hfill $\square$ \bigskip

\subsection{Proof of Theorem \ref{main}.}
\emph{Proof.} Denote $\og_{_N}^n=P_N(u(t_n))$ and $\eta^n=u_{_N}^n-\og_{_N}^n\in \widetilde{X}_N$. In view of \eqref{pro} and the triangle inequality, it is sufficient to show
\be\label{topr}
\|\eta^n\|\le M(\tau+h),
\ee
where $M$ is independent of $\tau$ and $h$ for $0\le n\tau\le T$.

The proof is given by induction. For $n=0$, it is obvious by using Lemma \ref{proint}:
\be\label{eta0}
\|\eta^0\|=\|I_N(u_0)-P_N(u_0)\|\le Ch\|u_0\|_1.
\ee
Suppose the claim is true for $n=0, 1, \ldots, k$. We prove that $\|\eta^{k+1}\|\le M(\tau+h)$. Subtracting \eqref{sch} from the projection of \eqref{local} in $X_N$ and noticing that the operator $P_N$ commutes with $e^{-\tau\p_x^3}$, we get for $n=0, 1,\ldots,k,$
\begin{align}
\eta^{n+1}&=e^{-\tau\p_x^3}\eta^n-\fl{\tau}{2} e^{-\tau\p_x^3}\left[I_N\delta_x^0 ((u_{_N}^n)^2)-P_N\delta_x^0(u(t_n)^2)\right]+\fl{c\tau h}{2}e^{-\tau\p_x^3}\delta_x^2 \eta^n-P_N(\xi^{n+1})\nn\\
&=e^{-\tau\p_x^3}\left[\eta^n+\fl{c\tau h}{2}\delta_x^2 \eta^n-\fl{\tau}{2}I_N\delta_x^0 \left((u_{_N}^n)^2-(\og_{_N}^n)^2\right)+\zeta^{n+1}\right],\label{eror}
\end{align}
where $$\zeta^{n+1}=\fl{\tau}{2}\left[P_N\delta_x^0(u(t_n)^2)-
I_N\delta_x^0((\og_{_N}^n)^2)\right]-e^{\tau\p_x^3}P_N(\xi^{n+1}).$$
It follows from Lemma \ref{proint}, \eqref{inp} and Lemma \ref{localer} that
\begin{align}
\|\zeta^{n+1}\|&\le \fl{\tau}{2}\left\|P_N\delta_x^0(u(t_n)^2)-I_N\delta_x^0(u(t_n)^2)\right\|+
\fl{\tau}{2}\left\|I_N\delta_x^0(u(t_n)^2-(\og_{_N}^n)^2)\right\|
+\|P_N(\xi^{n+1})\|\nn\\
&\le C\tau h\|\delta_x^0(u(t_n)^2)\|_1+\frac{\tau}{2}\interleave\delta_x^0\left(u(t_n)^2-
(\og_{_N}^n)^2\right)\interleave_{_N}+M_3(\tau^2+\tau h)\nn\\
&\le C\tau h\|u(t_n)^2\|_2+\frac{\tau}{2}\left\|u(t_n)^2-(\og_{_N}^n)^2\right\|_1+M_3(\tau^2+\tau h)\nn\\
&\le C \tau h\|u(t_n)\|_2^2+C\tau\|u(t_n)+\og_{_N}^n\|_1\|u(t_n)-\og_{_N}^n\|_1
+M_3(\tau^2+\tau h)\nn\\
&\le C \tau h\|u(t_n)\|_2^2+C\tau h\|u(t_n)\|_1\|u(t_n)\|_2
+M_3(\tau^2+\tau h)\le M_3(\tau^2+\tau h),\label{zeta}
\end{align}
where $M_3$ depends on $\|u\|_{L^\infty(0,T; H_{\rm p}^3(\Omega))}$.
Here for the third inequality we used the properties
\begin{align*}
\|\delta_x^0v\|^2&=\fl{1}{4h^2}\int_{-\pi}^{\pi}\left(v(x+h)-v(x-h)\right)^2dx
=\fl{1}{4h^2}\int_{-\pi}^{\pi}\left(\int_{-h}^hv'(x+y)dy\right)^2dx\\
&\le \fl{1}{2h}\int_{-\pi}^{\pi}\int_{-h}^h\left(v'(x+y)\right)^2dydx=|v|^2_1,\\
\interleave\delta_x^0 v\interleave^2_{_N}&=\fl{1}{4h}\sum\limits_{j=0}^{2N}\left(v(x_j+h)-v(x_j-h)\right)^2
=\fl{1}{4h}\sum\limits_{j=0}^{2N}\left(\int_{-h}^{h}v'(x_j+y)dy\right)^2\\
&\le \fl{1}{2}\sum\limits_{j=0}^{2N}\int_{-h}^{h}\left(v'(x_j+y)\right)^2dy=|v|_1^2,
\end{align*}
and the well-known bilinear estimate $\|fg\|_1\le C\|f\|_1 \|g\|_1$.
For simplicity of notation, we denote $u_j^n=u_{_N}^n(x_j)$, $\og_j^n=\og_{_N}^n(x_j)$, $\eta_j^n=\eta^n(x_j)$ and $\zeta_j^{n+1}=\zeta^{n+1}(x_j)$. Recall that $u_j^n$, $\og_j^n$, $\eta_j^n$, $\zeta_j^{n+1}\in \mathbb{R}$ by definition. Applying \eqref{eror}, Young's inequality, \eqref{inp} and \eqref{intp}, we obtain
\begin{align}
&\hspace{-3mm}\|\eta^{n+1}\|^2=\big\|\eta^n+\fl{c\tau h}{2}\delta_x^2 \eta^n-\fl{\tau}{2}I_N\delta_x^0 ((u_{_N}^n)^2-(\og_{_N}^n)^2)+\zeta^{n+1}\big\|^2\nn\\
&=\left\|\eta^n+\zeta^{n+1}\right\|^2+\fl{\tau^2}{4}\left\|ch\delta_x^2 \eta^n-I_N\delta_x^0 ((\eta^n)^2)-2I_N\delta_x^0(\eta^n\og_{_N}^n)\right\|^2\nn\\
&\quad+\tau\left(\eta^n+\zeta^{n+1}, ch\delta_x^2 \eta^n-I_N\delta_x^0 ((\eta^n)^2)-2I_N\delta_x^0(\eta^n\og_{_N}^n)\right)\nn\\
&\le (1+\tau)\|\eta^n\|^2+(1+\fl{1}{\tau})\|\zeta^{n+1}\|^2+\fl{\tau^2}{4}\interleave ch\delta_x^2 \eta^n-\delta_x^0 ((\eta^n)^2)-2\delta_x^0(\eta^n\og_{_N}^n)\interleave_{_N}^2\nn\\
&\quad+\tau\left\langle\eta^n+\zeta^{n+1}, ch\delta_x^2 \eta^n-\delta_x^0 ((\eta^n)^2)-2\delta_x^0(\eta^n\og_{_N}^n)\right\rangle_{_N}\nn\\
&=(1+\tau)\|\eta^n\|^2+(1+\fl{1}{\tau})\|\zeta^{n+1}\|^2+\fl{c^2\tau^2h^2}{4}
\interleave\delta_x^2\eta^n\interleave_{_N}^2
+\fl{\tau^2}{4}\interleave\delta_x^0((\eta^n)^2)\interleave_{_N}^2\nn\\
&\quad+\tau^2\interleave\delta_x^0(\eta^n\og_{_N}^n)\interleave_{_N}^2
-\fl{c\tau^2 h}{2}\left\langle\delta_x^2\eta^n, \delta_x^0((\eta^n)^2)\right\rangle_{_N}-c\tau^2 h\left\langle\delta_x^2\eta^n, \delta_x^0(\eta^n\og_{_N}^n)\right\rangle_N\nn\\
&\quad+\tau^2 \left\langle\delta_x^0((\eta^n)^2), \delta_x^0(\eta^n\og_{_N}^n)\right\rangle_N
-c\tau h\interleave\delta_x^+\eta^n\interleave_{_N}^2+c\tau h\left\langle\zeta^{n+1}, \delta_x^2 \eta^n\right\rangle_{_N}-\tau\left\langle\eta^n, \delta_x^0 ((\eta^n)^2)\right\rangle_{_N}\nn\\
&\quad-\tau\left\langle\zeta^{n+1}, \delta_x^0 ((\eta^n)^2)\right\rangle_{_N}-2\tau\left\langle\eta^n,\delta_x^0(\eta^n\og_{_N}^n)\right\rangle_{_N}
-2\tau\left\langle\zeta^{n+1}, \delta_x^0(\eta^n\og_{_N}^n)\right\rangle_{_N}.\label{sum}
\end{align}
Next we estimate the terms in \eqref{sum} separately by using similar arguments as in \cite{Courtes2017}. By definition and \eqref{eq0}, we have
\be\label{sp1}
\fl{c^2\tau^2h^2}{4}\interleave\delta_x^2 \eta^n\interleave_{_N}^2=\fl{c^2\tau^2h^3}{4}\sum\limits_{j=0}^{2N}\left(\delta_x^2\eta_j^n\right)^2\\
=c^2\tau^2\left(\interleave\delta_x^+\eta^n\interleave_{_N}^2-\interleave\delta_x^0
\eta^n\interleave_{_N}^2\right).
\ee
Moreover, it follows from \eqref{NI}, by induction $\|\eta^n\|\le M(\tau+h)$ and the assumption $\tau \le h/c$ that
\be\label{etaf}
\|\eta^n\|_{L^\infty}\le h^{-1/2}\|\eta^n\|\le M h^{-1/2}(\tau+h)\le M(1+1/c)h^{1/2}\le c,
\ee
whenever
\be\label{h1}
h\le h_1=M^{-2}c^4(1+c)^{-2}.
\ee
Thus, when $h\le h_1$ we have
\be
\fl{\tau^2}{4}\interleave\delta_x^0 ((\eta^n)^2)\interleave_{_N}^2=\fl{\tau^2 h}{4}\sum\limits_{j=0}^{2N}(\delta_x^0 \eta_j^n)^2(\eta_{j+1}^n+\eta_{j-1}^n)^2\le c\tau^2 \|\eta^n\|_{L^\infty} \interleave\delta_x^0 \eta^n\interleave_{_N}^2.\label{sp2}
\ee
In view of the Sobolev inequality and \eqref{pro}, we have
\be\label{ogp}
\|\og_{_N}^n\|_{L^\infty}\le C\|\og_{_N}^n\|_1\le C\|u(t_n)\|_1\le M_1,\quad
\|\p_x\og_{_N}^n\|_{L^\infty}\le C\|\og_{_N}^n\|_2\le C\|u(t_n)\|_2\le M_2,
\ee
where $M_1$ and $M_2$ depend on $\|u\|_{L^\infty(0,T;H_{\rm p}^1(\Omega))}$ and $\|u\|_{L^\infty(0,T; H_{\rm p}^2(\Omega))}$, respectively.
This yields
\begin{align}
\tau^2\interleave\delta_x^0(\eta^n\og_{_N}^n)\interleave_{_N}^2&=\tau^2 h\sum\limits_{j=0}^{2N}\left(\eta_{j+1}^n \delta_x^0\og_j^n+\og^n_{j-1}\delta_x^0 \eta^n_j\right)^2\nn\\
&\le 2\tau^2 h\sum\limits_{j=0}^{2N}\left((\eta_{j+1}^n)^2(\delta_x^0\og_j^n)^2+(\og^n_{j-1})^2(\delta_x^0 \eta_j^n)^2\right)\nn\\
&\le 2\tau^2\|\p_x\og_{_N}^n\|_{L^\infty}^2\|\eta^n\|^2+2\tau^2\|\og_{_N}^n\|^2_{L^\infty}
\interleave\delta_x^0\eta^n\interleave_{_N}^2\nn\\
&\le 2\tau^2\left(M_2^2\|\eta^n\|^2+M_1^2\interleave\delta_x^0\eta^n\interleave_{_N}^2\right).\label{sp3}
\end{align}
Applying \eqref{eq1} and \eqref{etaf}, we obtain
\begin{align}
-\fl{c\tau^2 h}{2}\left\langle\delta_x^2\eta^n, \delta_x^0((\eta^n)^2)\right\rangle_{_N}&=
\fl{c\tau^2}{2}\sum\limits_{j=0}^{2N}\eta_j^n\eta_{j+1}^n\delta_x^+\eta_j^n
-\fl{c\tau^2}{2}\sum\limits_{j=0}^{2N}\eta_{j-1}^n\eta_{j+1}^n\delta_x^0\eta_j^n\nn\\
&=-\fl{c\tau^2 h^2}{6}\sum\limits_{j=0}^{2N}(\delta_x^+\eta_j^n)^3+\fl{2c\tau^2h^2}{3}\sum\limits_{j=0}^{2N}
(\delta_x^0 \eta_j^n)^3\nn\\
&\le -\fl{c\tau^2 h^2}{6}\sum\limits_{j=0}^{2N}(\delta_x^+\eta_j^n)^3+\fl{2c\tau^2}{3}\|\eta^n\|_{L^\infty}
\interleave\delta_x^0\eta^n\interleave_{_N}^2.\label{sp4}
\end{align}
Similarly, using \eqref{eq2}, \eqref{ogp} and the assumption $c\tau\le h$ yields
\begin{align}
-c\tau^2 h\left\langle\delta_x^2\eta^n, \delta_x^0(\eta^n\og_{_N}^n)\right\rangle_N&=
c\tau^2\sum\limits_{j=0}^{2N}\eta_j^n\eta_{j+1}^n\delta_x^+\og_j^n-c\tau^2\sum\limits_{j=0}^{2N}
\eta_{j-1}^n\eta_{j+1}^n\delta_x^0 \og_j^n\nn\\
&\le 2\tau\|\p_x\og_{_N}^n\|_{L^\infty}\|\eta^n\|^2\le2 M_2 \tau\|\eta^n\|^2.\label{sp5}
\end{align}
Some tedious calculations give
\begin{align}
&\hspace{-5mm}\tau^2 \left\langle\delta_x^0((\eta^n)^2), \delta_x^0(\eta^n\og_{_N}^n)\right\rangle_N=\tau^2 h \sum\limits_{j=0}^{2N}(\delta_x^0\eta_j^n)(\eta_{j+1}^n+\eta_{j-1}^n)(\og_{j+1}^n\delta_x^0\eta_j^n
+\eta_{j-1}^n\delta_x^0\og_j^n)\nn\\
&=\tau^2h\sum\limits_{j=0}^{2N}(\delta_x^0\eta_j^n)^2\og_{j+1}^n(\eta_{j+1}^n+\eta_{j-1}^n)
+\tau^2 h\sum\limits_{j=0}^{2N}\eta_{j-1}^n(\eta_{j+1}^n+\eta_{j-1}^n)(\delta_x^0\eta_j^n) (\delta_x^0\og_j^n)\nn\\
&\le 2\tau^2\|\og_{_N}^n\|_{L^\infty}\|\eta^n\|_{L^\infty}\interleave\delta_x^0\eta^n\interleave_{_N}^2+
\fl{\tau^2}{h}\|\eta^n\|_{L^\infty}\|\p_x\og_{_N}^n\|_{L^\infty} h\sum\limits_{j=0}^{2N}|\eta_{j-1}^n|(|\eta_{j+1}^n|+|\eta_{j-1}^n|)\nn\\
&\le 2\tau^2\|\og_{_N}^n\|_{L^\infty}\|\eta^n\|_{L^\infty}\interleave\delta_x^0\eta^n\interleave_{_N}^2+2\tau
\fl{\tau}{h}\|\eta^n\|_{L^\infty}\|\p_x\og_{_N}^n\|_{L^\infty}\|\eta^n\|^2\nn\\
&\le 2 \tau^2M_1\|\eta^n\|_{L^\infty}\interleave\delta_x^0\eta^n\interleave_{_N}^2+
2M_2\tau\|\eta^n\|^2.\label{sp6}
\end{align}
Applying Young's inequality, we have
\begin{align}
&\hspace{-5mm}c\tau h\left\langle\zeta^{n+1}, \delta_x^2 \eta^n\right\rangle_{_N}=c\tau \sum\limits_{j=0}^{2N}\zeta_j^{n+1}\left(\eta_{j+1}^n-2\eta_j^n+\eta_{j-1}^n\right)\nn\\
&\le 2c\sum\limits_{j=0}^{2N}(\zeta_j^{n+1})^2+\fl{c\tau^2}{8}\sum\limits_{j=0}^{2N}
\left(\eta_{j+1}^n-2\eta_j^n+\eta_{j-1}^n\right)^2\nn\\
&\le 2c\sum\limits_{j=0}^{2N}(\zeta_j^{n+1})^2+2c\tau^2\sum\limits_{j=0}^{2N}
(\eta_{j}^n)^2\le\fl{2}{\tau}\|\zeta^{n+1}\|^2+2\tau\|\eta^n\|^2.\label{sp9}
\end{align}
Furthermore, a straightforward calculation yields that
\begin{align}
&\hspace{-5mm}-\tau\langle\eta^n, \delta_x^0((\eta^n)^2)\rangle_{_N}=-\fl{\tau}{2}\sum\limits_{j=0}^{2N} \eta_j^n\left((\eta_{j+1}^n)^2-(\eta_{j-1}^n)^2\right)=-\fl{\tau}{2}\sum\limits_{j=0}^{2N} \left[\eta_j^n(\eta_{j+1}^n)^2-\eta_{j+1}^n(\eta_{j}^n)^2\right]\nn\\
&=\fl{\tau}{6}\sum\limits_{j=0}^{2N} \left[(\eta_{j+1}^n)^3-3\eta_j^n(\eta_{j+1}^n)^2+3\eta_{j+1}^n(\eta_{j}^n)^2-(\eta_j^n)^3\right]
=\fl{\tau h^3}{6}\sum\limits_{j=0}^{2N}(\delta_x^+\eta_j^n)^3,\label{sp7}\\
&\hspace{-5mm}-2\tau\langle\eta^n, \delta_x^0(\eta^n\og_{_N}^n)\rangle_{_N}
= -\tau h\sum\limits_{j=0}^{2N}\eta_j^n\eta_{j+1}^n \delta_x^+\og_j^n\le \tau \|\p_x\og_{_N}^n\|_{L^\infty} \|\eta^n\|^2\le M_2\tau \|\eta^n\|^2.\label{sp8}
\end{align}
Similarly, one derives that
\begin{align}
&\hspace{-4mm}-\tau\left\langle\zeta^{n+1}, \delta_x^0 ((\eta^n)^2)\right\rangle_{_N}-2\tau\left\langle\zeta^{n+1}, \delta_x^0(\eta^n\og_{_N}^n)\right\rangle_{_N}\nn\\
&=-\fl{\tau}{2}\sum
\limits_{j=0}^{2N}\zeta_j^{n+1}\left[(\eta_{j+1}^n)^2-(\eta_{j-1}^n)^2+2\eta_{j+1}^n\og_{j+1}^n-
2\eta_{j-1}^n\og_{j-1}^n\right]\nn\\
&\le 2 \sum\limits_{j=0}^{2N}(\zeta_j^{n+1})^2+\fl{\tau^2}{4}\|\eta^n\|_{L^\infty}^2\sum\limits_{j=0}^{2N}
(\eta_{j+1}^n-\eta_{j-1}^n)^2+\fl{\tau^2}{4}\sum\limits_{j=0}^{2N}
\left(\eta_{j+1}^n\og_{j+1}^n-
\eta_{j-1}^n\og_{j-1}^n\right)^2\nn\\
&\le \fl{\tau}{h}\Big[\fl{2}{\tau} \|\zeta^{n+1}\|^2+\tau(\|\eta^n\|_{L^\infty}^2+\|\og_{_N}^n\|_{L^\infty}^2)\|\eta^n\|^2\Big]
\le \fl{1}{c}\Big[\fl{2}{\tau} \|\zeta^{n+1}\|^2+\tau(c^2+M_1^2)\|\eta^n\|^2\Big].\label{sp10}
\end{align}
Combining \eqref{sum} and \eqref{sp1}-\eqref{sp10}, we obtain that
\be\label{sum2}
\begin{split}
\|\eta^{n+1}\|^2&\le (1+A\tau)\|\eta^n\|^2+\left(1+\fl{3}{\tau}+\fl{2}{c\tau}\right)\|\zeta^{n+1}\|^2+B\tau^2\|\delta_x^0\eta^n\|^2\\
&\quad+\tau h\sum\limits_{j=0}^{2N}(h-c\tau)\left(\fl{h}{6}\delta_x^+\eta_j^n-c\right)(\delta_x^+\eta_j^n)^2,
\end{split}
\ee
where
\begin{align*}
A&=3+c+5M_2+2\tau M_2^2+M_1^2/c,\\
B&=2M_1^2-c^2+\|\eta^n\|_{L^\infty}\left(2M_1+5c/3\right).
\end{align*}
Applying \eqref{etaf}, we get
\[B\le 2M_1^2-c^2+2M(1+1/c)\left(M_1+c\right)h^{1/2},\]
which implies that $B\le 0$ whenever
\be\label{c0}
c> c_0=\sqrt{2}M_1,
\ee
and
\be\label{h0}
h\le h_0=\left(\fl{c^2-2M_1^2}{2M(1+1/c)(M_1+c)}\right)^2,
\ee
where $M_1$ is given by \eqref{ogp} depending on $\|u\|_{L^\infty(0,T;H_{\rm p}^1(\Omega))}$.
It is easily observed that $h_0\le h_1$.
In view of \eqref{etaf}, we have
\[\fl{h}{6}\delta_x^+\eta_j^n-c\le \fl{1}{3}\|\eta^n\|_{L^\infty}-c\le -\fl{2c}{3}<0,\quad \mathrm{if}\quad h\le h_0.\]
This together with the CFL condition $c\tau \le h$ and \eqref{sum2} yields that for $n=0,\ldots,k$,
\begin{align*}
\|\eta^{n+1}\|^2&\le (1+\tau C(M_2,c))\|\eta^n\|^2+\fl{C(c)}{\tau}\|\zeta^{n+1}\|^2\\
&\le (1+\tau C(M_2,c))\|\eta^n\|^2+\tau C(M_3,c)(\tau+h)^2,
\end{align*}
where $C(c,d)$ indicates that $C$ depends on $c$ and $d$.
Hence
\begin{align*}
\|\eta^{k+1}\|^2&\le e^{\tau C(M_2,c)}\|\eta^k\|^2+\tau C(M_3,c)(\tau+h)^2\\
&\le e^{2\tau C(M_2,c)}\|\eta^{k-1}\|^2+\tau C(M_3,c)(\tau+h)^2\big(1+
e^{\tau C(M_2,c)}\big)\\
&\le \ldots\\
&\le e^{(k+1)\tau C(M_2,c)}\|\eta^0\|^2+\tau C(M_3,c)(\tau+h)^2\big(1+
e^{\tau C(M_2,c)}+\ldots+e^{k\tau C(M_2,c)}\big)\\
&\le e^{(k+1)\tau C(M_2,c)}\big[\|\eta^0\|^2+\fl{C(M_3,c)}{C(M_2,c)}
(\tau+h)^2\big]\\
&\le e^{TC(M_2,c)}C(M_3,c)
(\tau+h)^2,
\end{align*}
which gives the error \eqref{topr} for $n=k+1$ by setting
\be\label{Md}
M=C(T, M_3, c)= e^{TC(M_2,c)/2}C^{1/2}(M_3,c).
\ee
This concludes the proof.

\section{Numerical experiments}
In this section, we present some numerical experiments to illustrate our analytic convergence rate given in Theorem \ref{main}. In practical computation, the interpolation $I_N$ is implemented via FFT, which is very efficient.

\bigskip

{\sl Example 1.} The well-known solitary-wave solution of the KdV equation \eqref{KdV} is given by
\be\label{sol-ex}
u(x,t)=12\lambda\,\sech^2(\sqrt{\lambda}(x-4\lambda t-a)),\quad a\in\mathbb{R},\quad \lambda>0.
\ee
It represents a single bump moving to the right with speed $4\lambda$. Here we choose $\lambda=1/4$ and the torus $\Omega=(-30,30)$ which is large enough such that the periodic boundary conditions do not introduce significant errors, i.e., the soliton is far enough away from the boundary for the considered time interval.

\begin{figure}[h!]
\begin{minipage}[t]{0.5\linewidth}
\centering
\includegraphics[width=3.15in,height=2.6in]{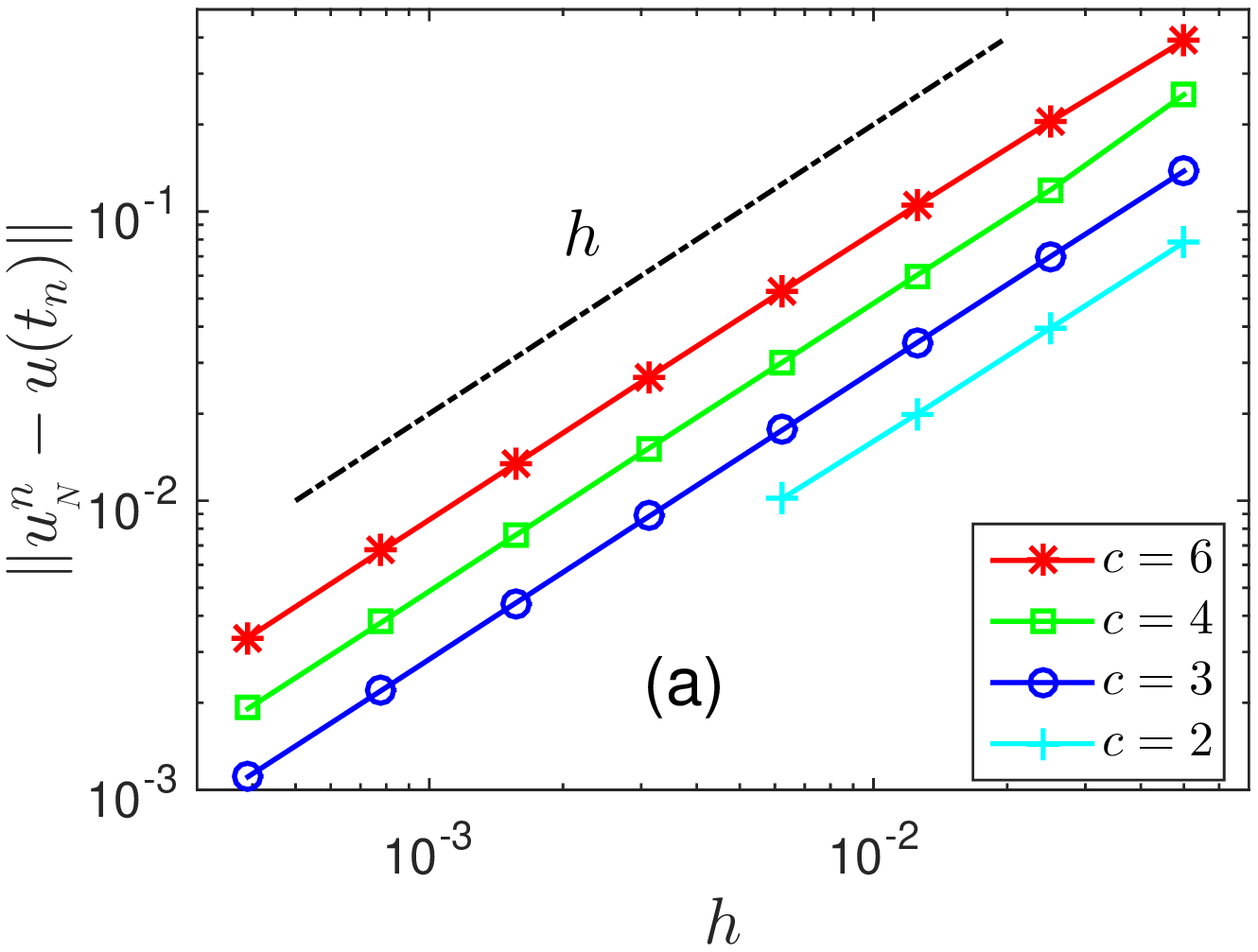}
\end{minipage}%
\hspace{3mm}
\begin{minipage}[t]{0.5\linewidth}
\centering
\includegraphics[width=3.15in,height=2.6in]{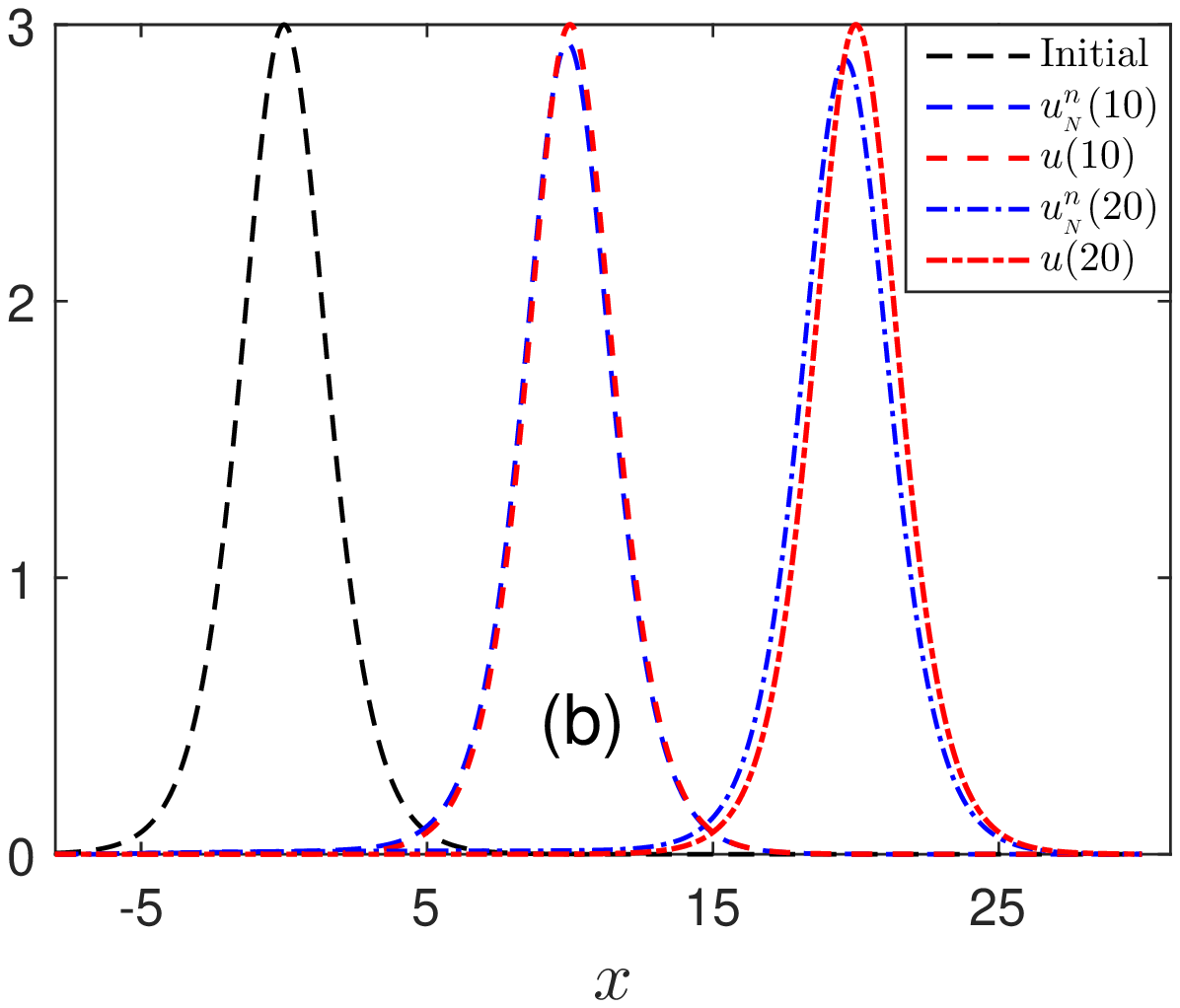}
\end{minipage}
\caption{Numerical simulation for the solitary-wave solution \eqref{sol-ex}. (a) The error of the first-order scheme \eqref{sch} at $T=2$ for various choices of $h$ and $c$. The time step size $\tau$ satisfies $\tau=h/c$. The broken line has slope one. (b) The numerical solution at $T=10, 20$ was obtained by the scheme \eqref{sch} with $h=1/200$ and $\tau=h/4$.}\label{fig1}
\end{figure}

Figure \ref{fig1}(a) displays the discretization errors for the scheme \eqref{sch} at $T=2$ for various choices of $h$ and $c$ with $\tau=h/c$. The results for $\tau=dh$ with $d\le 1/c$ are similar, which are omitted here for brevity. It can be clearly observed that the scheme \eqref{sch} converges linearly in space under the condition $\tau\le h/c$. Moreover, the error decreases as $c$ gets smaller, which is reasonable due to the fact that $c$ is the coefficient of the added artificial viscosity. The constraint of $c\ge c_0$ is verified by the fact that the numerical solution blows up when $h\le 1/320$ for $c=2$.
On the other hand, the solution also explodes when $\tau=dh$ with $d>1/c$, which shows the CFL condition $\tau\le h/c$ in Theorem \ref{main} is sharp. Figure \ref{fig1}(b) illustrates the time evolution of the solitary wave and the corresponding first-order approximate solutions for fixed $h=1/200$ and $\tau=h/4$.

\bigskip

\begin{figure}[h!]
\begin{minipage}[t]{0.5\linewidth}
\centering
\includegraphics[width=3.15in,height=2.6in]{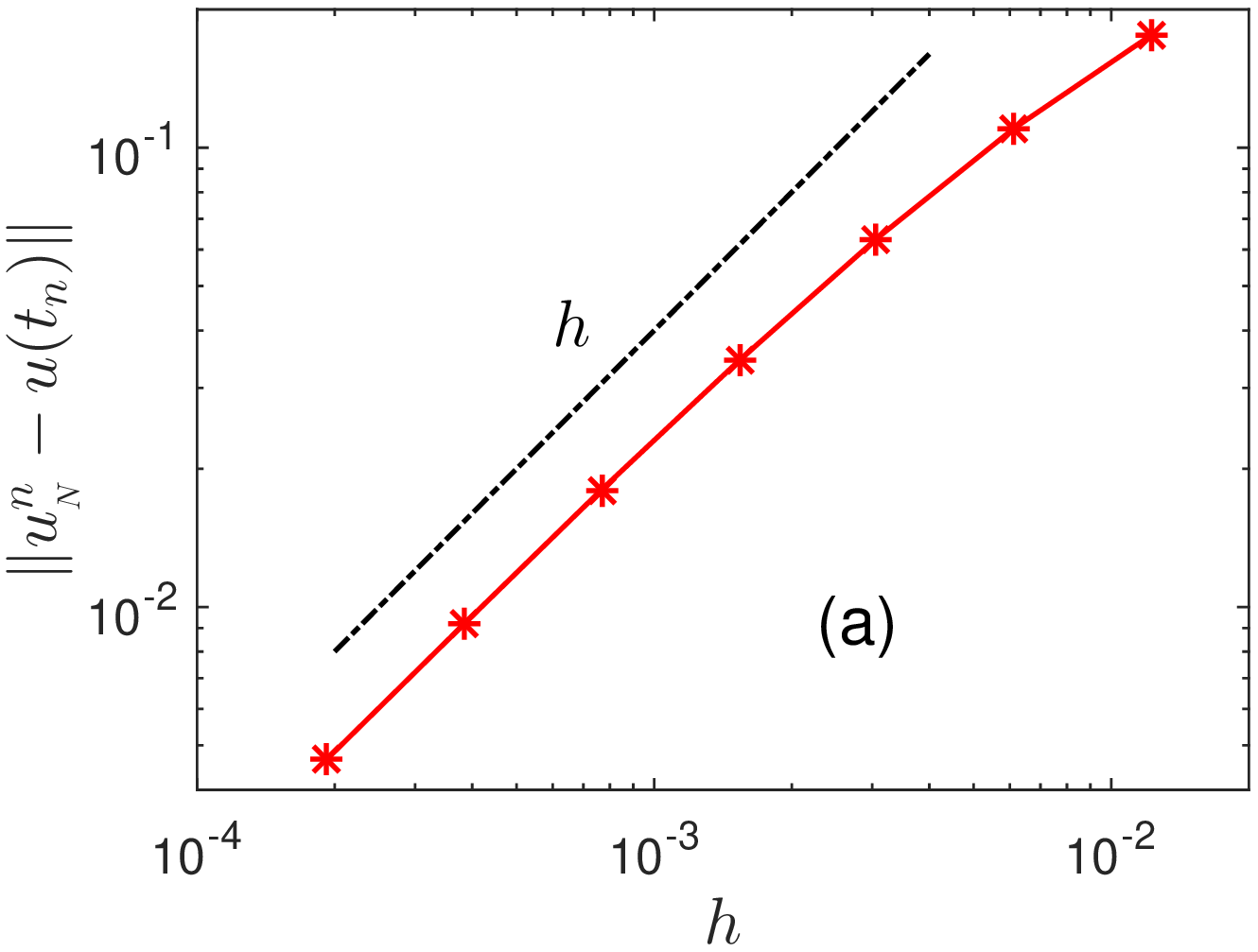}
\end{minipage}%
\hspace{3mm}
\begin{minipage}[t]{0.5\linewidth}
\centering
\includegraphics[width=3.15in,height=2.6in]{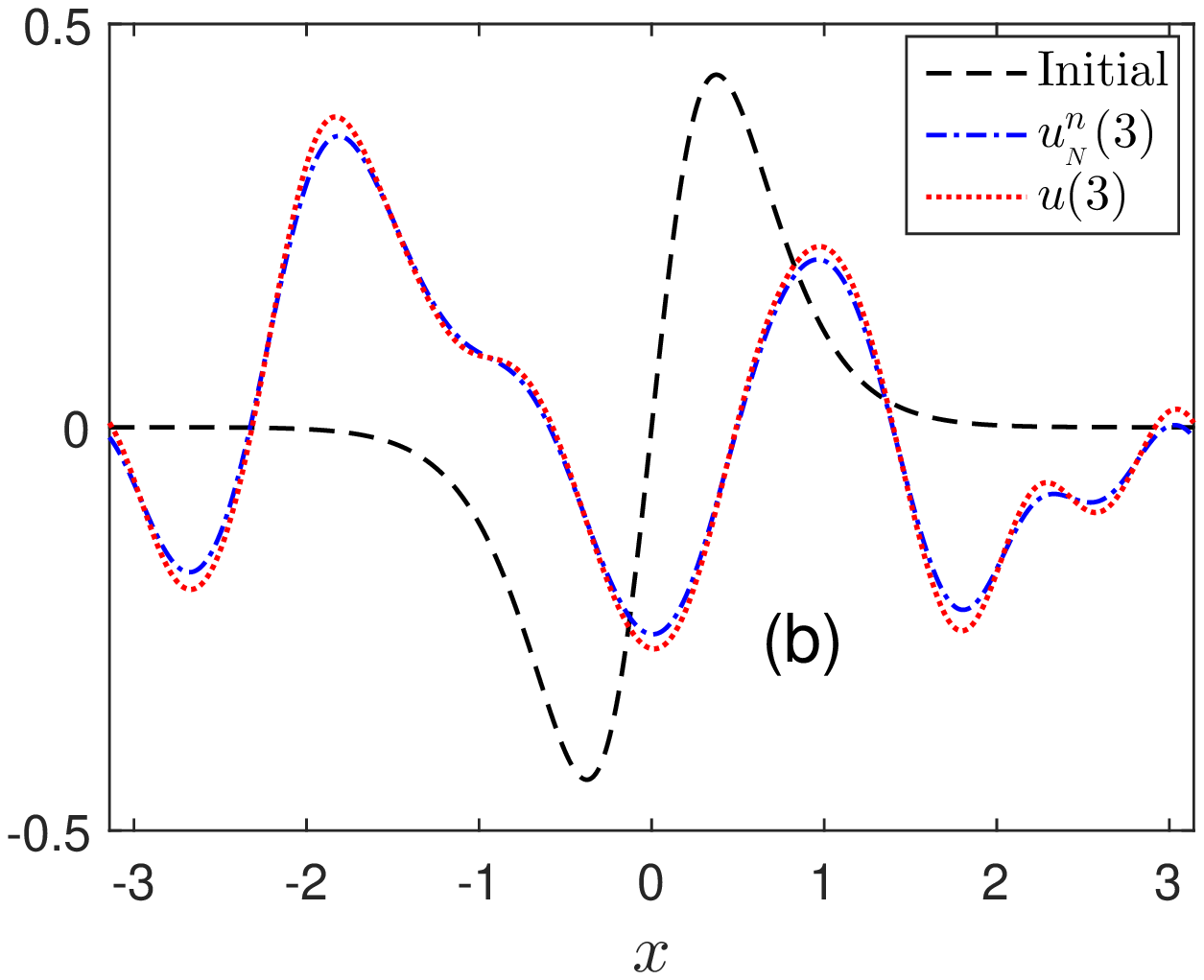}
\end{minipage}
\caption{Numerical simulation for the initial value \eqref{ex2}. (a) The error of the first-order scheme \eqref{sch} at $T=3$ for various choices of $h$ with $c=3$ and $\tau=h/\pi$. The broken line has slope one. (b) The numerical solution at $T=3$ was computed with the scheme \eqref{sch} using $h=\pi/2^{11}$ and $\tau=h/\pi$.}\label{fig2}
\end{figure}

{\sl Example 2.} The initial data of the KdV equation \eqref{KdV} is now chosen as
\be\label{ex2}
u_0(x)=3\,\sech^2(2x)\sin(x),\qquad x\in [-\pi, \pi].
\ee
The initial data and the numerical solution for $T=3$ with $c=3$, $h=\pi/2^{11}$ and $\tau=h/\pi$ are displayed in Figure \ref{fig2} (b), where the reference solution is obtained by the second-order exponential integrator of \cite{HS2017} with $\tau=10^{-6}$ and $h=\pi/2^{15}$.
The error of the scheme \eqref{sch} with $c=3$ and $\tau=h/\pi$ is shown in Figure \ref{fig2} (a). The graph clearly shows first-order convergence of the scheme \eqref{sch}.

\end{document}